\documentclass[12pt]{article}
\usepackage{amsmath,amsfonts,amssymb,amsthm}

\def\pf{\noindent{\bf Proof.\ \,}}
\def\eop{{$\square$}}
\def\labtt#1{\label {#1} }

\def\l{\lambda}
\def\o{\omega}

\def\QQ{{\mathbb Q}}

\def\ZZ{{\mathbb Z}}

\def\kron#1#2{\delta_{#1#2}}  

\def\la{\langle}
\def\ra{\rangle}
\def\vac{\hbox{\bf 1}} 

\def\vnat{V_{\natural}}

\begin{document}

\newtheorem{thm}{Theorem}[section]
\newtheorem{prop}[thm]{Proposition}
\newtheorem{lem}[thm]{Lemma}
\newtheorem{rem}[thm]{Remark}
\newtheorem{coro}[thm]{Corollary}
\newtheorem{conj}[thm]{Conjecture}
\newtheorem{de}[thm]{Definition}
\newtheorem{hyp}[thm]{Hypothesis}

\newtheorem{nota}[thm]{Notation}
\newtheorem{ex}[thm]{Example}
\newtheorem{proc}[thm]{Procedure}

\def\refpp#1{(\ref {#1})}

\def\bifo#1#2{{\langle #1, #2 \rangle}}

\begin{center}
{\Large \bf  Lattice-integrality of certain group-invariant integral forms in vertex operator algebras}

\end{center}

\centerline{ 18 May 2014  }
\bigskip

\begin{center}

Chongying Dong

Department of Mathematics,

University of California,

Santa Cruz, CA 95064 USA

{\tt dong@ucsc.edu}

\vskip 0.5cm
and
\vskip 0.5cm

Robert L. Griess Jr.

Department of Mathematics,

University of
Michigan,

Ann Arbor, MI 48109-1043  USA

{\tt rlg@umich.edu}

\smallskip

\end{center}

\begin{abstract}
Certain vertex operator algebras have integral forms (integral spans of bases which are closed under the countable set of products).    It is unclear when they (or integral multiples of them) are integral as lattices under the natural
bilinear form on the VOA.   We show that lattice-integrality may be arranged under some hypotheses, including cases of integral forms invariant by finite groups.    In particular, there exists a
lattice-integral Monster-invariant integral form in the Moonshine VOA.
\end{abstract}

\vskip 1cm


\newpage

\tableofcontents

\section{Introduction}

This article partially resolves an open question raised in \cite{ivoa} about lattice-integrality of an integral form in a vertex operator algebra invariant under a finite
group.  See the next subsection for definitions.   We note that integral valued form is obtained as an integer multiple of a given bilinear form which takes rational values  with bounded denominators.   The growth of denominators is not a trivial point in the infinite dimensional situation.

We present two general results for constructing IVOAs which are lattice-integral, \refpp{dongl1} and \refpp{quasiprimary}.
Also, we show that there is a Monster-invariant LIIVOA in the Moonshine module.    See \refpp{sharingef} and
\refpp{mvoaliif}.  Finally, we discuss some open questions about whether lattice-integrality can be arranged for integral forms in VOAs.

The earlier article \cite{ivoa} was concerned with integral forms invariant under a finite group.  In the present article, finite group actions play no role until
\refpp{hyp1}.

\subsection{Terminology}

Throughout this article, we use the following notations and definitions.   They are similar to usage in \cite{ivoa}.

Since the adjective ``integral'' has both a ring-theoretic meaning (as in integral form) and an arithmetic meaning (as in integral lattice), our notations are a bit heavy.

The dual of an additive subgroup $J$ in an inner product space $U$ is $J^{\circ}=\{u\in U \mid \bifo {u}{J}\subset \ZZ\}$.

The symbols $V, W$ stand for VOAs over the rationals.

An additive subgroup $A$ of $V$ is {\it lattice integral (LI)} if $\bifo AA \le \ZZ$, where $\bifo {\ }{\ }$ is the natural bilinear form \cite{fhl}  on $V$.  As usual, $\bifo {\vac}{\vac}=1$, where $\vac$ is the vacuum element of $V$.      We say that $A$ is  {\it nearly lattice integral} (NLI)  if there is an integer $m>0$ so
that $mA$ is lattice integral.

An {\it integral form} (IF) in $V$ shall mean an integral form of the vector space $V$ which is closed under all the VOA products and which
contains a positive integral multiple of the vacuum $\vac$ and a positive integer multiple of the principal Virasoro element.

An {\it integral VOA} (IVOA) shall mean an integral form  which contains the vacuum.    An IVOA is not assumed to be lattice-integral.

An IF,   $J$, which is lattice-integral is called  is an  LIIF.   If $J$ is nearly lattice integral, we say $J$ is a NLIIF. An IVOA which is LI, NLI is termed LIIVOA, NLIIVOA, respectively.

{\it Total eigenlattice} is defined in \refpp{tel}.

\section{Main Results}

We begin with a few results about creating  integral forms and adjusting integral forms to ones with nicer properties.
\subsection{From NLIIF to LIIF}

\begin{nota} For subsets $X, Y$ of the VOA $V$, we define $X\cdot Y$ to be the set of all $a_kb$, for all $a \in X, b\in Y, k \in \ZZ$.
\end{nota}

\begin{lem}\labtt{1}  Suppose that $V$ is a VOA with a nondegenerate invariant bilinear form,  $V_0=\QQ\vac$, $(\vac, \vac )=1$
and  $J\subset V$ is a homogeneous IVOA. If $m>0$ is an integer so that $mJ$ is LI, then $m  J +\ZZ \vac$ is a LIIVOA.
\end{lem}
\pf
We first  show that $K:=m J +\ZZ \vac$ is an IVOA. For $a,b\in J$ and $n\in\ZZ$ then $(ma)_nmb=m^2a_nb\in m^2J\subset mJ$ and
$(ma)_n\vac=m(a_n\vac)\in m J.$   Finally, $\vac_n K \le K$ \refpp{voaprod1}.  Therefore, $K$ is an  IVOA.

Now we show that inner products are integral.   We have $\bifo {mJ}{mJ} \le \ZZ$ by assumption.
To compute $\bifo K {\vac}$, observe that this is the same as the set of all $\bifo {Px}{\vac}$, where $P:V\rightarrow \QQ \vac$ is the
orthogonal projection.
Homogeneity of $J$ implies that $P(J)= J \cap \QQ \vac \le \ZZ \vac$ \refpp{zzvac}.   Now use the hypothesis that $\bifo \vac \vac=1.$
\eop
\begin{lem} \labtt{dualstable} Suppose that $V=\oplus_{n\geq 0} V_n$ is a VOA with $\dim  (V_0)=1$ and that $V$ has a nondegenerate invariant bilinear form.
Assume that  $J\subset V$ is a homogeneous IVOA. Then $\frac{L(1)^n}{n!}J^{\circ}\subset J^{\circ}$ for all
$n\geq 0$
where $J^{\circ}=\{u\in V|\bifo {u}{J}\subset \ZZ\}$ is the dual of $J$ in $V.$
\end{lem}
\pf
Note that for any $u\in J$ and $n\geq 0,$  $\frac{L(-1)^n}{n!}u=u_{-n-1}\vac \in J$ as $\vac\in J$ (see \refpp{voaprod1}(ii)).    So $\frac{L(-1)^n}{n!}$ is well defined on $J.$  Then for any $u\in J^{\circ}, v\in J$ and $n\geq 0$
$$\bifo{\frac{L(1)^n}{n!}u}{v}=\bifo{u}{\frac{L(-1)^n}{n!}v}\in\ZZ.$$
That is, $\frac{L(1)^n}{n!}u\in J^{\circ},$ as desired.
\eop

\begin{lem}\label{dongl1} Suppose that $V=\oplus_{n\geq 0} V_n$ is a VOA with $\dim  (V_0)=1$ and that $V$ has a nondegenerate invariant bilinear form. Assume that  $J\subset V$ is an IVOA generated by $\sum_{s\leq t}J_s.$
Then there exists $m>0$ such that $\ZZ\vac+ m\sum_{s=1}^tJ_s$ generates a LIIVOA.
\end{lem}
\pf Clearly, $J^{\circ}$ is also homogeneous. Notice that $\sum_{s\leq
t}J_s\cap J^{\circ}_s$ has finite index in
both $\sum_{s\leq t}J_s$ and $\sum_{s\leq t}J^{\circ}_s.$ Let $m_1,m_2>0$ such that
$$ m_1\sum_{s\leq t}J_s\subset \sum_{s\leq t}J_s\cap J^{\circ}_s, \ \ \ \ m_2\sum_{s\leq t}J^{\circ}_s\subset \sum_{s\leq t}J_s\cap J^{\circ}_s.$$
Then for any $u\in J_a$ ($a\leq t$) we have, by \refpp{dualstable},  $m_1m_2\frac{L(1)^n}{n!}u \in \sum_{s\leq t}J_s\cap J^{\circ}_s.$

Let $m=m_1m_2.$ Then,
$\frac{L(1)^n}{n!}m\sum_{s\leq t}J_s\subset \sum_{s\leq t}J_s\cap J^{\circ}_s\subset \sum_{s\leq t}J_s$
and $\frac{L(1)^n}{n!}m\sum_{s\leq t}J_s\subset \sum_{s\leq t}J^{\circ}_s$ for all $n\geq 0.$ It is clear that
$\bifo{u}{v}\in\ZZ$ for $u\in J_0+m\sum_{s=1}^tJ_s, v\in \sum_{s\leq t}J_s.$ Note that $J_0=\ZZ\vac$ \refpp{zzvac}.

Denote by $J(m)$ the IF generated by $J_0+m\sum_{s=1}^tJ_s.$ Since $J$ is an IVOA by hypothesis, the sub-IF $J(m)$ is an IVOA.

Note that $J(m)=\sum_{s}J(m)_s.$ We shall prove by induction on $s\ge 0$ that
$$\bifo{J(m)_s}{J(m)_s}\in\ZZ.$$
For $s=0$, we use $J_0=\ZZ\vac$ and $\la \vac, \vac \ra = 1$.

We observe that since $J$ is an IVOA, for all $k$, we have $J_k(mJ) \le mJ$.  Consequently, $J(m)_\ell \le mJ_\ell$ for all $\ell$ and 
$J(m)_\ell = mJ_\ell$ for $1 \le \ell \le t$.    
Therefore,   for $s\in \{1,2, \dots, t\}$,  we have $\bifo{J(m)_s}{J(m)_s} \le m\ZZ$, whence 
$$\bifo{J(m)_s}{J(m)_s}\in\ZZ \leqno (*)$$ for all $s, 0\le s \le t$.

Now assume that $s>0$ and that (*) holds for all nonnegative indices up to $s$.   We now prove the result for index $s+1$.   Let $u,v\in J(m)_{s+1}.$ It is enough
to prove that $\bifo{u}{v}\in\ZZ$ in case $u$ has shape
$u=a_{-n-1}\vac$ or, for some $k\ge 1$, shape $u=a^1_{n_1}\cdots a^{k}_{n_k}a$ for homogeneous $a,a^1,...,a^k\in  m\sum_{s\leq t}J_s$ and
 $n_i\in\ZZ$.

Let $u=a_{-n-1}\vac.$ Then by \refpp{dualstable},
 $$\bifo{a_{-n-1}\vac}{v}=\bifo{\frac{L(-1)^n}{n!}a}{v}=\bifo{a}{\frac{L(1)^n}{n!}v}\in\ZZ.$$
For $u=a^1_{n_1}\cdots a^{k}_{n_k}a$, with $a$ homogeneous, we consider
\begin{eqnarray*}
& &\bifo{Y(a^1,z_1)\cdots Y(a^k,z_k)a}{v}\\
& &\ \ \ \ =\bifo{a}{Y(e^{L(1)z_k}(-z_k^{-2})^{L(0)}a^k,z_k^{-1})\cdots Y(e^{L(1)z_1}(-z_1^{-2})^{L(0)}a^1,z_1^{-1})v}
\end{eqnarray*}

Since the right argument is in $\ZZ[[z_1,z_1^{-1},\cdots z_k,z_k^{-1}]]J(m)$ and since $a\in J(m)$ is homogeneous of degree $d \le s$, the value of this inner product lies in
 $\ZZ[[z_1,z_1^{-1},\cdots z_k,z_k^{-1}]] \la J(m)_d, J(m)_d \ra$, which, by induction, is contained in  $\ZZ[[z_1,z_1^{-1},\cdots z_k,z_k^{-1}]]$.   The proof is complete.
\eop

\subsection{IVOAs and quasi-primary vectors}

Recall that $v\in V$   is called a {\it quasi primary vector}  if $L(1)v=0.$  Let $QP(V)$ be the subspace
of $V$ consisting of quasi primary vectors of $V.$

\begin{lem}\label{quasiprimary}Assume that  $J\subset V$ is an IVOA generated by quasi primary vectors.
Then $J$ is a LIIVOA.
\end{lem}

\pf Let $J$ be generated by a homogeneous quasi primary vectors $P.$  Then $J$ is a $\ZZ$-span of
$a^1_{n_1}\cdots a^k_{n_k}\vac$ for $a^i\in P$ and $n_i\in\ZZ.$ Then for any $v\in J$
we have
\begin{eqnarray*}
& &\bifo{Y(a^1,z_1)\cdots Y(a^k,z_k)\vac}{v}\\
& &\ \ \ \ =\bifo{\vac}{Y(e^{L(1)z_k}(-z_k^{-2})^{L(0)}a^k,z_k^{-1})\cdots Y(e^{L(1)z_1}(-z_1^{-2})^{L(0)}a^1,z_1^{-1})v}\\
& &\ \ \ \ =\bifo{\vac}{Y((-z_k^{-2})^{L(0)}a^k,z_k^{-1})\cdots Y((-z_1^{-2})^{L(0)}a^1,z_1^{-1})v}.
\end{eqnarray*}
Since $Y((-z_k^{-2})^{L(0)}a^k,z_k^{-1})\cdots Y((-z_1^{-2})^{L(0)}a^1,z_1^{-1})v\in J[[z_1,z_1^{-1},\cdots z_k,z_k^{-1}]]$ and $\bifo{\vac}{\vac}=1$,   the result follows immediately from \refpp{zzvac}. 
\eop

Lemma \ref{quasiprimary} can be applied to many examples. Recall from \cite{ivoa} that the lattice vertex operator algebra $V_L$ associated to the positive definite even lattice $L=\oplus_{i=1}^d\ZZ\gamma_i$ has an IVOA $(V_L)_{\ZZ}$ generated by $e^{\pm\gamma_i}$ for $i=1,...,d.$ It was proved directly in \cite{ivoa} that $\bifo{(V_L)_{\ZZ}}{(V_L)_{\ZZ}}\leq \ZZ.$ Since $e^{\pm\gamma_i}$ for $i=1,...,d$ are primary vectors for the Virasoro algebra, applying Lemma  \ref{quasiprimary} to $(V_L)_{\ZZ}$ to reprove this result.

\subsection{Ubiquity of NLIIFs}

Here we show how to adjust a NLI to a more convenient one, then make applications to group-invariant IFs.   Existence of a Monster-invariant lattice-integral IVOA in the Moonshine VOA follows.   The existence theory for IVOAs in \cite{ivoa}  had left the lattice-integrality property generally unresolved.

\begin{lem}\labtt{2}   Let $W$ be a VOA over $\QQ$.  Suppose $K$ is a finitely generated IF in $W$.

(i) If $J$ is a finitely generated IF in $W$, there is an integer $m>0$ so that $mJ \le K$,

(ii)  If $K$ is NLI, so is $J$.

(iii)  Given two finitely generated IFs $J, J'$ in $W$,  the abelian groups  $J/J\cap J'$ and $J'/J\cap J'$ have finite exponent.
Furthermore, $J$ is NLI if and only if $J'$ is NLI.
\end{lem}
\pf
Let $K$ be a finitely generated IF in $W$.   We may assume that $K$ is homogeneous.   Choose a finite set of integers, $S$, so  that
$K_S:=\sum_{i \in S} K_i$ generates $K$ as an IF.

Now let $J$ be an IF in $W$.   Since $J$ is finitely generated, we may enlarge $S$ by a finite amount to get a set of generators for $J$ in
$W_S:=\sum_{i \in S} W_i$.    Since $W_S$ is a finite dimensional rational space, there exists an integer $m>0$ so that
$m(J\cap W_S)\le K$.  It follows that $mJ \le K$.
This proves (i) and (ii).

Now, we prove (iii).   There exist integers $m>0$ and $n>0$ so $mJ\le K$ and so
$J/(J\cap K)$ has exponent dividing $m$ and
$J'/(J'\cap K)$ has exponent dividing $n$.  It follows that
$(J+J')/(J\cap J'\cap K)$ has exponent dividing $mn$.   The equivalence of property NLI for $J, J'$  now follows.
\eop

Now, we introduce group actions to play a role as in \cite{ivoa}.

\begin{hyp}\labtt{hyp1}
Suppose that $V$ is a VOA over the rationals and that the finite group
$N\le Aut(V)$ contains a normal subgroup $E\cong 2^r$, for $r\ge 2$ and, by conjugation, $N$ induces acts transitively on the nontrivial linear
characters of $E$.
Take $D \le E$ and suppose that $rank(E/D)\ge 1$.   Define $W:=V^D$, fixed point subVOA.   Let $J$ be an integral form in $V$.
\end{hyp}

\begin{rem} There a similar but stronger hypothesis in \cite{ivoa}, Notation 5.3, requiring $r\ge 3$ and the action of $N$ to be doubly
transitive.   Part of that item is repeated in \refpp{hyp2}, below.
\end{rem}\begin{thm}\labtt{sharingef}  We use the notations \refpp{hyp1}.   Suppose that $W$ has a finitely generated IF $K$ which is NLI.

(i) If $J$ is an IF in $V$ then $J$ contains a finitely generated NLIIF, $J'$.

(ii) Furthermore, if $G$ is a finite subgroup of $Aut(V)$ so that $J$ is $G$-invariant, then we may choose $J'$ to be $G$-invariant and
furthermore to satisfy
$J'\cap \QQ \vac = J\cap \QQ \vac$ and $J'\cap \QQ \o = J\cap \QQ \o$.
\end{thm}
\pf
(i) We may replace $J$ by the subVA generated over $\ZZ$ by
$\sum_{i \in S}J_i$, where $S$ is a finite set of indices.  Invariance under $G$ is preserved.

To prove that $J$ is $NLI$, we may without loss assume that $J$ is $N$-invariant \refpp{inv1}.   Then $J^D$ (fixed points) is an $IF$ in $W$.
By \refpp{2}, there is an integer $m>0$ so that
the IF $mJ^D$ is LI.

If $\l$ is a complex linear character of $E$, it follows that whenever $D \le Ker(\l)$, $mJ_{\l}$  is LI.   Since $N$ induces the action of
$GL(r,2)$ on $E$, this is true for every character $\l$ of $E$.   Therefore $mTel(J)$ is a IF and is LI.

Finally we observe that $2^r J \le Tel(J):=\sum_{\l \in \hat E} J_l$.   It follows that the IF $m2^r J$ is LI.
\eop

\begin{hyp} \labtt{hyp2}
Recall that
$V$ is a VOA over the rationals.

A group
$N\le Aut(V)$ has a normal subgroup $E\cong 2^r$ for $r\ge 3$;
 $N$ acts doubly transitively on the nontrivial linear characters of $E$;
 define   $N_{\l} :=Stab_N(\l )$.

We have $D$, a  subgroup of $E$ and $rank(E/D)\ge 2$;
 there exists an IVOA  $A$ in the fixed point subVOA $W:=V^D$
  and $A$ is invariant by $N_N(D)$
\end{hyp}

Under hypothesis \refpp{hyp2} and the assumption that $dim(V_0)=1$, we showed in \cite{ivoa} that  there exists an $N$-invariant IVOA $B$
in $V$.    Furthermore, if we have a finite group $G\le Aut(V)$ so that $N\le G$, then $B$ may be chosen to be $G$-invariant.

\begin{thm}\labtt{maintheorem}
We use the notations of \refpp{hyp2} and $B$ as above.  Suppose that $V$ is finitely generated.    If $A$ is NLI, then we may arrange for $B$ to be NLI.
\end{thm}
\pf  If necessary, replace $B$ by a finitely generated subIVOA of $B$,
then use \refpp{2}(ii).
\eop

The standard IVOA $IV_L$ in a VOA of lattice type is LI.   Any finitely generated IVOA in $V_L$ is NLI.

\begin{coro}\labtt{mvoaliif}
The Monster leaves invariant a LIIVOA in the Moonshine VOA.
\end{coro}
\pf  It suffices to prove that there is a homogeneous invariant NLIIVOA \refpp{1}.

As in the proof of Theorem 5.9 of \cite{ivoa},
we use the standard IVOA  in $(\vnat )^z \cong V_{\Lambda}^+$ (where $z$ is a $2B$ involution and $\Lambda$ is the Leech lattice).  This IVOA
is LI.    It follows that the Monster-invariant IVOA created in Theorem 5.9 of \cite{ivoa} is NLI.
\eop

\section{Open Questions}

\hskip 0.6 cm
{\bf OQ1:} If $J$ is an IF in a VOA, is there an integer $m>0$ so that $(J, J) \le \frac 1m \ZZ$?  If yes, then $mJ$ is an IF which is also
lattice-integral.   (It may help to assume that $J$ is finitely generated.)

{\bf OQ2:}   Assume that the IF $J$ is homogeneous, $J=\bigoplus_i J_i$.    For $a, b \in J_i$, $a_{i-1}b \in J_i$.   On $J_i$, each $a_{i-1}$
is represented by an integer matrix, so we have an integer valued form on $J$, the direct sum of the forms
$Trace( (a_{i-1}|_{J_i}) (b_{i-1}|_{J_i}) )$ on the $J_i$.
What is the relationship between this form and the rational valued form on $J$ inherited from $V$?   In general, these forms are not proportional and the second form is not generally associative.

For an example of the latter, look at the dihedral VOA of type $2A$, in degree $i=2$.   This degree 2 component has dimension 3 and is spanned by three conformal vectors of central charge $\frac 12$.   Details follow.
We use the following basis for the dihedral algebra of type $2A$ (this would be the degree 2 part of a VOA generated by a certain pair of conformal vectors of central charge 1/2).   The following comes from \cite{sakuma}, using formulas from  \cite{sym4}, page 2449.

Type 2A, basis $a,b,c$, with inner product and algebra product as follows:
For every permutation $p,q,r$ of $a,b,c$ we have
$(p,p)=1, (p,q)=\frac 18$ and
$p\cdot p=p$, $p\cdot q = \frac 18 (p+q-r)$.

It follows that $ad(a)$ has matrix  (for left action, using basis $a,b,c$)
$A:= \begin{pmatrix}
1 & \frac 18 & \frac 18 \cr
0 &  \frac 18 & -\frac 18 \cr
0 & -\frac 18 & \frac 18 \cr
\end{pmatrix}$
and
that $ad(b)$ has matrix
$B:= \begin{pmatrix}
\frac 18 & 0 & -\frac 18 \cr
\frac 18 & 1 &  \frac 18 \cr
-\frac 18 & 0 & \frac 18 \cr
\end{pmatrix}$.

We used Maple for the matrix work which follows.
We have $A^2=\begin{pmatrix}
1&\frac 18&\frac 18\cr
0&\frac 1{32}&-\frac 1{32}\cr
0&-\frac 1{32}&\frac 1{32}
\end{pmatrix}$
and
$AB=\begin{pmatrix}
\frac 18&\frac 18&-\frac 3{32}\cr
\frac 1{32}&\frac 18&0\cr
-\frac 1{32}&-\frac 18&0
\end{pmatrix}$.
We calculate the matrix for $ad(c)$, $C:=\begin{pmatrix}
\frac 18&-\frac 18 & 0\cr
-\frac 18 & \frac 18 & 0\cr
\frac 18 & \frac 18 & 1
\end{pmatrix}$ and $AC=\begin{pmatrix}
\frac 18 & -\frac 3{32} & 0\cr
-\frac 1{32} & 0 & -\frac 18 \cr
\frac 1{32} & 0 & \frac 18
\end{pmatrix}$.
The ``Killing form'' $\kappa (x,y):=Tr(ad(x)ad(y))$ therefore has Gram matrix
$\begin{pmatrix}
\frac{17}{16}&\frac 14&\frac 14\cr
\frac 14 & \frac{17}{16}& \frac 14\cr
\frac 14 & \frac 14 & \frac{17}{16}
\end{pmatrix}$.
The natural bilinear form, $\nu$,  (described earlier) has Gram matrix
$\begin{pmatrix}
1&\frac 18&\frac 18\cr
\frac 18 & 1& \frac 18\cr
\frac 18 & \frac 18 & 1
\end{pmatrix}$.
Since these two Gram matrices are not proportional, neither are $\kappa$ and $\nu$.

Another property to consider is associativity of a bilinear form.   The form $\nu$ is associative.   For $\kappa$, consider $ad(a^2)ad(b)$, which has matrix $AB$, and
$ad(a)ad(ab)$ which has matrix $A\cdot \frac 18(A+B-C) =
\frac 18
 \begin{pmatrix}
1 & \frac 18 & \frac 18 \cr
0 &  \frac 18 & -\frac 18 \cr
0 & -\frac 18 & \frac 18 \cr
\end{pmatrix}
\begin{pmatrix}
1 & \frac 14 & 0 \cr
\frac 14 & 1 & 0 \cr
-\frac 14 & -\frac 14 & -\frac 34
\end{pmatrix}$, which equals
$\frac 18 \begin{pmatrix}
1& \frac {11}{32} & -\frac 3{32} \cr
\frac 1{16} & \frac 5{32} &\frac 3{32} \cr
-\frac 1{16} & -\frac 5{32} & -\frac 3{32}
\end{pmatrix}$.
Since $AB$ has trace $\frac 14$ and $A\cdot \frac 18(A+B-C)$ has trace $\frac {17}{128}$, $\kappa$ is not associative.    This gives a second proof that $\kappa$ and $\nu$ are not proportional.

\appendix

\section{Appendix:   Background results}

\begin{de}\labtt{tel}  Given the action of an elementary abelian $2$-group $E$ on the free abelian group $L$, the eigenlattice $L_{\l}$ associated to
character $\l \in Hom(E, \{\pm 1\} )$, is
$x\in L \mid gx=\l (g)x, \text{ for all }g \in E \}$.   The total eigenlattice is $Tel(L):=\bigoplus_{\l \in Hom(E, \{\pm 1\} )} L_{\l}$.
\end{de}

A general reference for total eigenlattices associated to actions of finite groups on lattices is \cite{gal}.

\begin{lem}\labtt{exponentpowerof2}
$2^r$ annihilates the quotient group $J/Tel(J)$.
\end{lem}
\pf
Use the formula  $2^{-r}\sum_{g \in E} \l (g)g$ for the idempotent in the group algebra $\QQ [E]$ associated to $\l$.
\eop

\begin{lem}\labtt{inv1} (i)
If $G$ is a finite subgroup of $Aut(V)$ and $J$ is an IF, then $\bigcap _{g\in G} gJ$ is a $G$-invariant IF and $J/\bigcap _{g\in G} gJ$ is a
torsion abelian group.

(ii) If in addition, $J$ is finitely generated, there exists a $G$-invariant IF, $Q$, and an integer $n>0$ so that $nJ \le Q$.
\end{lem}
\pf
(i)
The invariance by $G$ is obvious.
Since $G$ is finite and each $gJ$ is an integral form of the rational vector space $V$,
$J/\bigcap _{g\in G} gJ$ is torsion.

(ii)
In case $J$ is finitely generated,
let $K$ be a finite rank additive subgroup of $J$ which generates $J$.

Then there exists an integer $m>0$ so that the finite abelian quotient group
$K/\bigcap _{g\in G} gK$ has exponent $m$.    Then $mJ \le subVA(mK) \le subVA(\bigcap _{g\in G} gK)$, and the latter is $G$-invariant.
\eop

\begin{lem}\labtt{zzvac}
If $J$ is a  IF, $J\cap \QQ\vac = n\ZZ \vac$ for some integer $n>0$.
\end{lem}
\pf
(A similar result  is in \cite{ivoa}, Remark 5.6.)
By definition of  IF, $J$ is a free abelian group and
$J\cap \QQ\vac$ contains a positive integer multiple of $\vac$,
so $J\cap \QQ\vac =  r \ZZ \vac$, for some rational number $r>0$.
Suppose that $r=a/b$ is not an integer, where $a$ and $b\ge 2$ are comprime integers.
The property  $\vac _{-1}\vac =\vac$ implies that $J$ contains $r^2 \vac$, a contradiction since $r \ZZ \vac$ is a direct summand of the free
abelian group $J$.
\eop

\begin{lem}\labtt{voaprod1} In a VOA $V$, let $a\in V$ and let $\vac$ be the vacuum.  Then for an integer $k \in \ZZ$,

\noindent
(i)
 $\vac_k=\kron{-1,}{k} 1_V$, i.e.  $\vac_k a=\kron{-1,}{k}a$ ;
(ii) $a_k \vac = \begin{cases} 0& k \ge 0\cr a & k=-1 \cr \frac 1{(k-1)!} (L_{-1})^{-k-1}a & k \le -2\cr \end{cases}$.
\end{lem}
\pf
See \cite{flm}.
\eop

\end{document}